\documentclass[conference]{ieeeconf}
\renewcommand{\baselinestretch}{0.998}

\renewcommand{\textfloatsep}{5mm} %slightly reduces space under figure captions

\IEEEoverridecommandlockouts
\usepackage[acronym]{glossaries}
\usepackage{amsmath,amssymb,amsfonts}
\usepackage{algorithmic}
\usepackage[pdftex]{graphicx}
\usepackage{xcolor}
\def\BibTeX{{\rm B\kern-.05em{\sc i\kern-.025em b}\kern-.08em
    T\kern-.1667em\lower.7ex\hbox{E}\kern-.125emX}}

\usepackage{units}
\usepackage{mathtools} 
\usepackage{amsthm}
\usepackage{lettrine}

\usepackage{scrextend} %required for list of abbreviations
\usepackage{tabularx}
\usepackage[english]{babel}
\usepackage{cite}

\usepackage{array}
\usepackage{color}
\usepackage{amssymb}
\usepackage{multicol}

\DeclareUnicodeCharacter{2212}{-}

\usepackage{todonotes}

\usepackage{tikz}
\usetikzlibrary{arrows,calc}
\usepackage{pgfplots}
\usepgfplotslibrary{groupplots}

\usepackage{amsthm}
\usepackage{algorithm}
\usepackage{hyperref}
\usepackage{epstopdf}

\newcounter{thm}
\newtheoremstyle{mystyle}%                % Name
  {}%                                     % Space above
  {}%                                     % Space below
  {}%                                     % Body font
  {}%                                     % Indent amount
  {\bfseries\color{black}}%               % Theorem head font
  {.}%                                    % Punctuation after theorem head
  {\newline}%                                    % Space after theorem head, ' ', or \newline
  {\thmname{#1}\thmnumber{ #2: }\normalfont\color{black}\thmnote{ \textit{(#3)}}}%                                     % Theorem head spec (can be left empty, meaning `normal')
\theoremstyle{mystyle}
\newtheorem{prob}[thm]{Problem}

\newtheorem{theorem}{Theorem}[section]

\newacronym{acr:cvt}{CVT}{continuously variable transmission}
\newacronym{acr:dp}{DP}{dynamic programming}
\newacronym{acr:ecms}{ECMS}{equivalent consumption minimization strategies}
\newacronym{acr:eltms}{ELTMS}{equivalent lap time minimization strategies}
\newacronym{acr:em}{EM}{electric motor}
\newacronym{acr:es2k}{ES2K}{Energy Storage to Kinetic}
\newacronym{acr:F1}{F1}{Formula 1}
\newacronym{acr:FIA}{FIA}{F\'{e}d\'{e}ration Internationale de l'Automobile}
\newacronym{acr:fgt}{FGT}{fixed-gear transmission}
\newacronym{acr:FD}{FD}{final drive}
\newacronym{acr:ice}{ICE}{internal combustion engine}
\newacronym{acr:k2es}{K2ES}{Kinetic to Energy Storage}
\newacronym{acr:mgu}{MGU}{motor generator unit}
\newacronym{acr:mguh}{MGU-H}{motor generator unit heat}
\newacronym{acr:mguk}{MGU-K}{motor generator unit kinetic}
\newacronym{acr:mpc}{MPC}{model predictive control}
\newacronym[description={Energy Management Strategy}, \glslongpluralkey={Energy Management Strategies},\glsshortpluralkey={EMSs}]{EMS}{EMS}{Energy Management Strategy}%
\newacronym{acr:pmp}{PMP}{Pontryagin's Minimum Principle}
\newacronym{acr:pu}{PU}{power unit}
\newacronym[description={Powertrain Operation}, \glslongpluralkey={Powertrain Operations},\glsshortpluralkey={POs}]{acr:PO}{PO}{Powertrain Operation}%
\newacronym{acr:rmse}{RMSE}{root-mean-square error}
\newacronym{acr:socp}{SOCP}{second-order cone program}
\newacronym{acr:soe}{SoE}{State of Energy}

\newcommand{\pushright}[1]{\ifmeasuring@#1\else\omit\hfill$\displaystyle#1$\fi\ignorespaces}
\newcommand{\pushleft}[1]{\ifmeasuring@#1\else\omit$\displaystyle#1$\hfill\fi\ignorespaces}
\makeatother
\newif\ifmargincomments %A quick way of turning off margin comments for, say, arXiv submission
\margincommentstrue
\newif\ifextendedversion %A quick way of turning off appendix
\extendedversionfalse
\ifmargincomments

\else

\fi
\begin{document}

\title{\bf A Sequential Convex Programming Approach to Free-trajectory Minimum-lap-time Optimization of Racing Cars}

%{Was a bit hesitant to claim "operation optimization" when only powertrain operation and racing line follow from optimizer, not steering angle. Also put free-trajectory before minimum-lap-time because free implies its being optimized already.}
	
	%Fast free-trajectory minimum lap-time optimization for racing cars


\author{Erik van den Eshof, Wytze de Vries, Jorn van Kampen, Mauro Salazar%
\thanks{%Erik van den Eshof, Jorn van Kampen, and Mauro Salazar are with the 
	The authors are with the Control Systems Technology section, Department of Mechanical Engineering, Eindhoven University of Technology (TU/e), Eindhoven, 5600 MB, The Netherlands.
E-mails: \tt\footnotesize r.c.p.v.d.eshof@tue.nl, w.a.b.d.vries@tue.nl,  j.h.e.v.kampen@tue.nl, m.r.u.salazar@tue.nl.}
}

\maketitle

\begin{abstract}
	
%\msmargin{The problem of simulating a racing car to drive as fast as possible over a race course is typically computationally expensive and, as a result, not feasible for real-time applications unless major simplifications are imposed. Among these simplifications is the assumption of a fixed path along the racing track, which forms a major limitation in achievable lap-time.}{Keep it to the point. What is the question?}
This paper presents a modeling and optimization framework to compute the minimum-lap-time spatial trajectory and powertrain operation of racing cars in a computationally efficient fashion. Specifically, we first derive a quasi-steady-state model of a racing car, whereby the racing line trajectory is jointly optimized.
Next, we frame the minimum-lap-time problem and leverage its mostly convex structure by devising a sequential convex programming solution algorithm.
We benchmark our method against off-the-shelf nonlinear programming solvers, showing how it can bring computation time down from a few minutes to a few seconds, paving the way for real-time implementations.
Moreover, we compare our results to similarly efficient minimum-curvature racing line optimization methods, showing how a minimum-time-based racing line might lead to 4\% faster lap-times.
Finally, we showcase our framework for optimal powertrain energy management and we validate the common modeling assumption that the racing line is unaffected by energy limitations, showing that this assumption results in marginal lap-time losses of under 0.1\%.
%Moreover, we use our framework to benchmark simpler models energy management models based on a maximum speed profile, validating their adequacy .
%  and to \msmargin{energy management models based on a maximum speed profile, showing how such an assumption results into a marginal lap-time loss.}{do we want to say it?}

%This paper presents a minimum lap-time optimization strategy that co-optimizes the driven trajectory in a computationally efficient fashion.
%\msmargin{Through a sequential-convex programming method on a specialized model, we achieve a minimum-time optimization and trajectory planning strategy with a solving time in the order of seconds, whereas traditional nonlinear-programming methods take over a minute to solve. Moreover, we demonstrate that a direct minimum-time-based trajectory yields significant lap-time benefits of several seconds over a similarly computationally efficient method based on a minimum-curvature objective.}{rewrite} %In addition to trajectory optimization, our method is capable of finding optimal powertrain operation, and we investigate a common assumption that energy management strategy is not significantly affected by driving trajectory.
\end{abstract}

\vspace{-0.5mm}
\section{Introduction}
\begin{figure}[t]
	\centering
	%\vspace{2mm} %so it aligns with abstract
	\includegraphics[width=\linewidth]{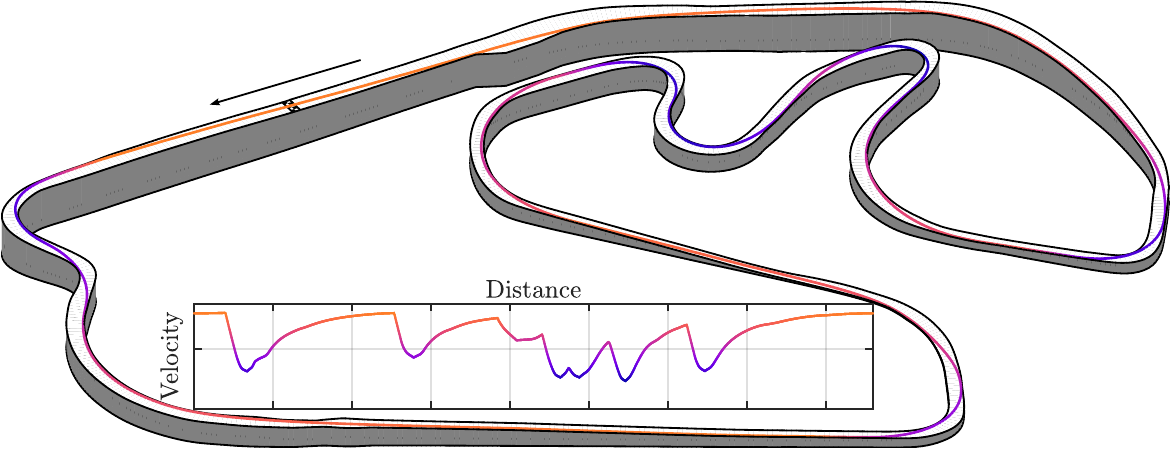}
	\caption{Optimal racing line and speed trajectory over the Interlagos Circuit.}
	\label{fig:track}
	\vspace{-2.2mm}
\end{figure}

%\lettrine{M}{otorsport} is an industry where every millisecond counts: in a relentlessly competitive environment, the identification of design solutions and operational strategies that maximize performance is of paramount importance.
\lettrine{M}{otorsport} is an industry of tiny margins: in this relentlessly competitive environment, maximizing vehicle performance is essential.
Identifying opportunities for improvement requires computational models that are capable of accurately and efficiently predicting a racing car's ability to cover the distance around a circuit, a \textit{lap}, in minimal time (Fig. \ref{fig:track}).
Typically, this problem poses significant computation time challenges, prompting existing approaches to simplify the models, sacrificing accuracy and ultimately achievable performance.
This paper addresses this challenge by proposing a modeling and optimization framework that retains key model considerations---most notably, the freedom to optimize the driving trajectory (the \textit{racing line})---whilst maintaining computational efficiency.
\paragraph*{Related Literature}
We can principally distinguish between two approaches to lap-time optimization: efficient \textit{fixed-trajectory} methods that assume a fixed racing line, and computationally heavy \textit{free-trajectory} methods where the racing line is jointly optimized~\cite{Massaro2021}.
In fixed-trajectory optimization, the racing line trajectory is computed in preprocessing~\cite{Braghin2008,Heilmeier2019} or based on onboard measurements~\cite{EbbesenSalazarEtAl2018}. These methods typically use quasi-steady-state (QSS) models. The simplest method is the so-called \textit{apex-finding} method, which finds the critical grip-limited speeds using a predetermined curvature trajectory and applies forward-backward integration to satisfy acceleration and deceleration limitations, resulting in a speed profile~\cite{Lovato2021,Massaro2021}. Alternatively, convex models are applied to efficiently capture more detailed vehicle dynamics~\cite{EbbesenSalazarEtAl2018,BorsboomFahdzyanaEtAl2021,Duhr2022,Kampen2024,EshofKampenEtAl2025}, or by directly exploiting vehicular data by constructing a (convex) \textit{g-g-v} performance envelope~\cite{Lovato2021,Lovato2021a,DuhrChristodoulouEtAl2020,Milliken1995}.
Reported solving times for these models are in the order of seconds. Yet these methods are not suitable to optimize the spatial trajectory due to the nonlinear coupling between curvature, velocity, and vehicular forces.

Alternatively, nonlinear programming (NLP) optimization models that jointly optimize the race trajectory often involve transient dynamics and fewer simplifications~\cite{Lovato2021a,Limebeer2018,Perantoni2014,Limebeer2015,Rucco2015,Limebeer2022,Gabiccini2021,Christ2019}. Although advances in computational capabilities have reduced the computation times of these models to the order of minutes, they remain unsuitable for real-time applications. 

In an attempt to find a middle ground between accuracy and efficiency, several sources propose a QSS model with a free trajectory~\cite{Lovato2021a,Veneri2019}, but reported solving times are still over a minute. Kapania et al. propose an iterative optimization scheme that alternates between minimizing curvature and minimizing lap-time~\cite{Kapania2019}, but the racing line in this method primarily determined by minimizing its curvature, which does not necessarily lead to the fastest lap-time~\cite{Heilmeier2019}.

In conclusion, to the best of the authors' knowledge, no existing method achieves an efficient joint optimization of the global driving trajectory and powertrain operation within a pure minimum-lap-time framework. 

\setlength{\fboxsep}{0pt}
\setlength{\fboxrule}{1pt}
\paragraph*{Statement of Contributions}
This paper presents an optimization model and solution framework based on sequential convex second-order conic programming to compute minimum-time driving trajectories and powertrain operational strategies for racing cars.
%Thereby we obtain solving times similar to fixed-trajectory methods with the accuracy benefits of minimum-time-based free-trajectory methods. Several nonlinear equations that could not be captured within the convex framework are linearized and subsequently solved through an iterative sequential-convex programming (SCP) approach.
As most of our model preserves convexity, the proposed sequential solution approach results in a fast convergence of a few seconds.
Additionally, we showcase the lap-time benefits of our minimum-time-based trajectory approach over an equally efficient minimum-curvature-based method, and demonstrate its potential for energy management applications.
\section{Methodology}
\label{ModelSection}
In this section, we illustrate the minimum-lap-time control problem and propose an iterative, sequential-convex programming solution method to retain a fast solving time. We build on an existing convex lap-time optimization framework commonly used for energy management~\cite{Duhr2022,EbbesenSalazarEtAl2018,Kampen2024,EshofKampenEtAl2025,MurgovskiJohannessonEtAl2015}, written as a second-order conic program (SOCP). As our focus is on vehicle dynamics, we extend it with load-dependent tire grip coefficients, lateral load transfer effects, and cornering resistances, while retaining convexity. We include a variable trajectory racing line, which introduces a number of nonconvex constraints. We convexify these constraints through linearizations around the result of a previous sequential-programming iteration. Thereby, we capture most of the complexity of the problem in the convex subproblem, with the rationale that this reduces the required number of subproblem evaluations for solving.\\
\subsection{Circuit Model}
\begin{figure}[t]
	\centering
	\includegraphics[width=0.75\linewidth]{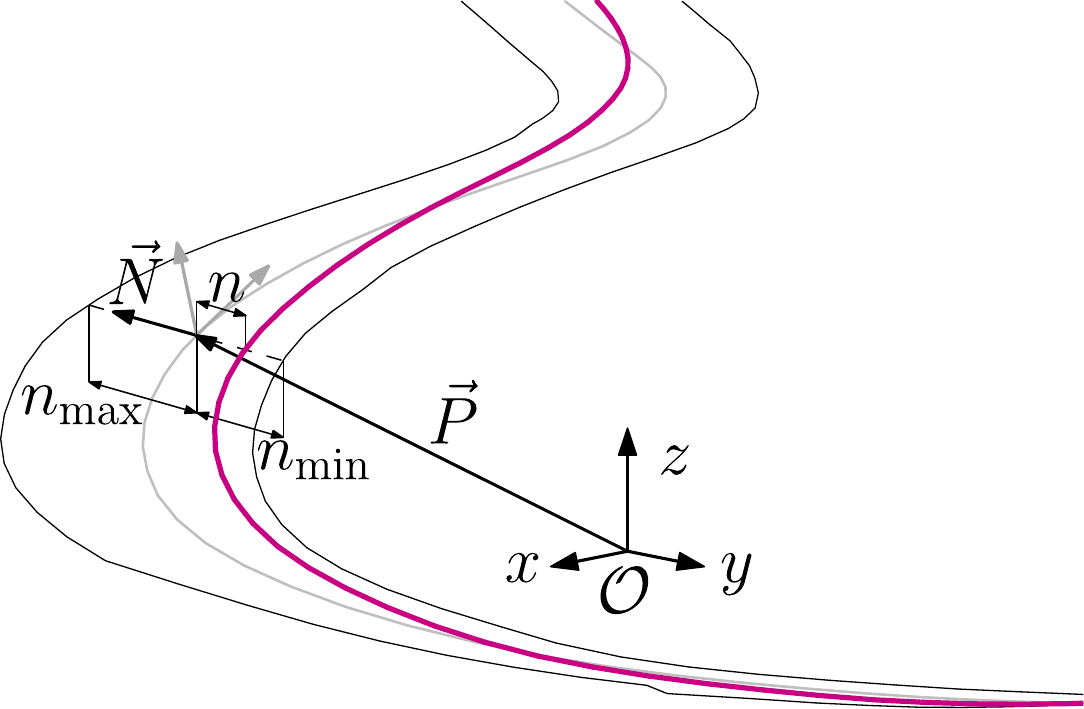}
	\caption{3D-ribbon track model, defined by a central trajectory and a normal vector. Positions on the track are defined by an offset along this vector.}
	\label{fig:trackmodel}
\end{figure}

To model the position and curvature of the car's trajectory along the circuit, we leverage the commonly used three-dimensional \textit{ribbon} model~\cite{Heilmeier2019,Perantoni2015}. As shown in Fig.~\ref{fig:trackmodel}, positions along the track are defined by a lateral offset of a track reference line, i.e,
\par\nobreak\vspace{-5pt}
\begingroup
\allowdisplaybreaks
\begin{small}

\begin{equation}
\label{eqn1}
\left[\begin{array}{c}
     x(s_\mathrm{ref})\\
     y(s_\mathrm{ref})\\
     z(s_\mathrm{ref})\\
\end{array}\right] = \underbrace{\left[\begin{array}{c}
     x_\mathrm{ref}(s_\mathrm{ref})\\
     y_\mathrm{ref}(s_\mathrm{ref})\\
     z_\mathrm{ref}(s_\mathrm{ref})\\ 
\end{array}\right]}_{\vec{P}(s_\mathrm{ref})} + n(s_\mathrm{ref})\underbrace{\left[\begin{array}{c}
     N_\mathrm{x}(s_\mathrm{ref})\\
     N_\mathrm{y}(s_\mathrm{ref})\\
     N_\mathrm{z}(s_\mathrm{ref})\\ 
\end{array}\right]}_{\vec{N}(s_\mathrm{ref})},
\end{equation}
\end{small}%
\endgroup
where $[x,y,z]^\top$ are the positional coordinates of the vehicle trajectory, $s_\mathrm{ref}$ is the distance along the circuit reference path, $\vec{P}=[x_\mathrm{ref},y_\mathrm{ref},z_\mathrm{ref}]^\top$ are the circuit reference path coordinates, $\vec{N}$ is the unit-length normal vector, and $n$ is the lateral offset along this vector. $\vec{P}$ and $\vec{N}$ are circuit-specific and subject to identification~\cite{Perantoni2015}, but fixed in the context of the optimization problem. The circuit boundaries are constrained on the lateral offset
\par\nobreak\vspace{-5pt}
\begingroup
\allowdisplaybreaks
\begin{small}

\begin{equation}
\label{eqn2}
n(s_\mathrm{ref}) \in \left[n_\mathrm{min}(s_\mathrm{ref})+\frac{w_\mathrm{veh}}{2},n_\mathrm{max}(s_\mathrm{ref})-\frac{w_\mathrm{veh}}{2}\right],
\end{equation}
\end{small}%
\endgroup
where $w_\mathrm{veh}$ is the width of the vehicle, to allow for sufficient space to the boundaries.\\
For notational simplicity, we denote the distance derivatives along the circuit reference with an apostrophe, i.e., $\frac{\mathrm{d}(\cdot)}{\mathrm{d}s_\mathrm{ref}} = (\cdot)'$. Accordingly, the first and second derivatives of the traveled trajectory follow from chain rule differentiation:
\par\nobreak\vspace{-5pt}
\begingroup
\allowdisplaybreaks
\begin{small}
\begin{equation}
\label{eqn3}
\left[\begin{array}{c}
     x'(s_\mathrm{ref})\\
     y'(s_\mathrm{ref})\\
     z'(s_\mathrm{ref})\\
\end{array}\right] = \vec{P}'(s_\mathrm{ref}) + n(s_\mathrm{ref})\vec{N}'(s_\mathrm{ref}) 
 + n'(s_\mathrm{ref})\vec{N}(s_\mathrm{ref})
\end{equation}
\end{small}%
\endgroup
and
\par\nobreak\vspace{-5pt}
\begingroup
\allowdisplaybreaks
\begin{small}
\begin{multline}
\label{eqn4}
\left[\begin{array}{c}
     x''(s_\mathrm{ref})\\
     y''(s_\mathrm{ref})\\
     z''(s_\mathrm{ref})\\
\end{array}\right] = \vec{P}''(s_\mathrm{ref}) + n(s_\mathrm{ref})\vec{N}''(s_\mathrm{ref}) 
 + n''(s_\mathrm{ref})\vec{N}(s_\mathrm{ref})\\ + 2n'(s_\mathrm{ref})\vec{N}'(s_\mathrm{ref}).
\end{multline}
\end{small}%
\endgroup
All of these equations are convex over the independent variable $n$ and its derivatives. However, to describe centrifugal force, we introduce a nonconvex relation.
\par\nobreak\vspace{-5pt}
\begingroup
\allowdisplaybreaks
\begin{small}
\begin{equation}
\label{eqn5}
F_\mathrm{c}(s_\mathrm{ref}) = 2 E_\mathrm{kin}(s_\mathrm{ref}) \underbrace{\frac{y''(s_\mathrm{ref})x'(s_\mathrm{ref})-x''(s_\mathrm{ref})y'(s_\mathrm{ref})}{\left(x'^2(s_\mathrm{ref})+y'^2(s_\mathrm{ref})\right)^{\frac{3}{2}}}}_{\kappa(s_\mathrm{ref})},
\end{equation}
\end{small}%
\endgroup
where $\kappa$ is the curvature and $E_\mathrm{kin}$ the kinetic energy of the vehicle. For simplicity's sake, we neglect the second-order torsional effects and vertical curvature by considering a planar curvature in the $xy$-plane for our centrifugal force formulation, which is a reasonable assumption as slope and banking changes remain small on typical racing circuits. This allows us to consider the important three-dimensional effects of slope and banking~\cite{Lovato2021,Lovato2021a}, while keeping the problem well-conditioned. The slope of the trajectory is defined as follows, again nonconvex:
\par\nobreak\vspace{-5pt}
\begingroup
\allowdisplaybreaks
\begin{small}
	\begin{equation}
		\label{eqn7nl}
		\theta(s_\mathrm{ref}) = \arctan{\left(\frac{z'(s_\mathrm{ref})}{\sqrt{x'^2(s_\mathrm{ref})+y'^2(s_\mathrm{ref})}}\right)}.
	\end{equation}
\end{small}%
\endgroup
The banking angle is defined as
\par\nobreak\vspace{-5pt}
\begingroup
\allowdisplaybreaks
\begin{small}
\begin{equation}
\label{eqn8}
\phi(s_\mathrm{ref}) = \arctan{\left(\frac{N'_{\mathrm{z}}(s_\mathrm{ref})}{\sqrt{N'^2_{\mathrm{x}}(s_\mathrm{ref})+N'^2_\mathrm{y}(s_\mathrm{ref})}}\right)}.
\end{equation}
\end{small}%
\endgroup

\subsection{Objective}
As our framework disregards potential competitor interactions, the objective is simply to traverse the circuit as quickly as possible and minimize lap-time. We are operating in the circuit reference distance domain, thereby our nonconvex objective is
\par\nobreak\vspace{-5pt}
\begingroup
\allowdisplaybreaks
\begin{small}
\begin{equation}
\label{eqn9}
    \mathrm{min}\, t_\mathrm{lap} = \mathrm{min} \int_{0}^{S_\mathrm{ref}}{\textstyle\frac{\mathrm{d}t}{\mathrm{d}s}}(s_\mathrm{ref}){\textstyle\frac{\mathrm{d}s}{\mathrm{d}s_\mathrm{ref}}}(s_\mathrm{ref})\, \mathrm{d}s_\mathrm{ref},
\end{equation}
\end{small}%
\endgroup
where $t_\mathrm{lap}$ is the total time required to complete the lap, $S_\mathrm{ref}$ is the length of the circuit, ${\textstyle\frac{\mathrm{d}t}{\mathrm{d}s}}$ is the \textit{lethargy}, the inverse of the vehicle velocity, and $s$ without the subscript is the distance traveled by the vehicle.  ${\textstyle\frac{\mathrm{d}s}{\mathrm{d}s_\mathrm{ref}}}$ follows from the circuit positional coordinates through the Pythagorean theorem:
\par\nobreak\vspace{-5pt}
\begingroup
\allowdisplaybreaks
\begin{small}
\begin{equation}
\label{eqn11}
     {\textstyle\frac{\mathrm{d}s}{\mathrm{d}s_\mathrm{ref}}} (s_\mathrm{ref})\geq \left\|\begin{array}{c}
x'(s_\mathrm{ref})\\
y'(s_\mathrm{ref})\\
z'(s_\mathrm{ref})\\
\end{array}\right\|_2.
\end{equation}
\end{small}%
\endgroup
Note that we have relaxed this equation to ensure convexity. Due to the presence of $\frac{\mathrm{d}s}{\mathrm{d}s_\mathrm{ref}}$ in the objective, this inequality holds with equality in an optimal solution. Furthermore, to model the inverse relation between lethargy and vehicle velocity $v$, we use the following relaxed constraint:
\par\nobreak\vspace{-5pt}
\begingroup
\allowdisplaybreaks
\begin{small}
\begin{equation}
\label{eqn12}
    {\textstyle\frac{\mathrm{d}t}{\mathrm{d}s}}(s_\mathrm{ref}) + v(s_\mathrm{ref}) \geq \left\|\begin{array}{c}
2\\
{\textstyle\frac{\mathrm{d}t}{\mathrm{d}s}}(s_\mathrm{ref}) - v(s_\mathrm{ref})
\end{array}\right\|_2,
\end{equation}
\end{small}%
\endgroup
which corresponds to 
\par\nobreak\vspace{-5pt}
\begingroup
\allowdisplaybreaks
\begin{small}
	\begin{equation}
		\left(\textstyle\frac{\mathrm{d}t}{\mathrm{d}s}(s_\mathrm{ref}) \geq \frac{1}{v(s_\mathrm{ref})}\right) \wedge \left(v(s_\mathrm{ref}) \geq 0\right)
	\end{equation}
\end{small}%
\endgroup
and holds with equality for an optimal solution as the lethargy is minimized~\cite{EbbesenSalazarEtAl2018,MurgovskiJohannessonEtAl2015}. Finally, we obtain the relation for kinetic energy through the vehicle velocity
\par\nobreak\vspace{-5pt}
\begingroup
\allowdisplaybreaks
\begin{small}
\begin{equation}
\label{eqn13}
     E_\mathrm{kin}(s_\mathrm{ref}) \geq \frac{1}{2}mv^2(s_\mathrm{ref}),
\end{equation}
\end{small}%
\endgroup
where $m$ is the vehicle mass. This again holds with equality for an optimal solution, as $v$ is maximized.
\subsection{Vehicle Dynamic Model}
\begin{figure}[t]
    \centering
    \includegraphics[width=\linewidth]{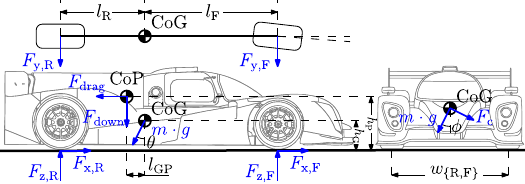}
    \caption{Top, side and frontal view of the simulated model. The forces acting on the vehicle are highlighted. The arrows indicate positive force directions.}
    \label{fig:1}
\end{figure}
The vehicle-dynamic model is visualized in Fig.~\ref{fig:1}. We employ a single-track model with steady-state cornering and dynamics~\cite{Pacejka2012,Milliken1995}. The coordinate frame used for forces is in the vehicle reference frame. In the remainder of this paper, we completely omit the distance dependency ($s_\mathrm{ref}$) from our notations. The longitudinal force balance is given as
\par\nobreak\vspace{-5pt}
\begingroup
\allowdisplaybreaks
\begin{small}
\begin{equation}
\label{eqn14}
{\textstyle\frac{\mathrm{d}E_\mathrm{kin}}{\mathrm{d}s}} = F_\mathrm{x,R} + F_\mathrm{x,F} - \cos(\phi)\,\sin(\theta)\,m\,g - \underbrace{\frac{1}{m}\,\rho \,C_\mathrm{d}A\,E_\mathrm{kin}}_{F_\mathrm{drag}},
\end{equation}
\end{small}%
\endgroup
where $F_\mathrm{x,R}$ and $F_\mathrm{x,F}$ are the longitudinal forces on the rear and front axles, respectively, $\rho$ is the air density, $C_\mathrm{d}$ is the drag coefficient, $A$ is the cross-sectional aerodynamic reference area, and $g$ is the gravitational acceleration. To map this derivative to our reference domain, we require a nonconvex relation
\par\nobreak\vspace{-5pt}
\begingroup
\allowdisplaybreaks
\begin{small}
\begin{equation}
\label{eqn15nl}
{\textstyle\frac{\mathrm{d}E_\mathrm{kin}}{\mathrm{d}s_\mathrm{ref}}} = {\textstyle\frac{\mathrm{d}E_\mathrm{kin}}{\mathrm{d}s}}\cdot {\textstyle\frac{\mathrm{d}s}{\mathrm{d}s_\mathrm{ref}}}.
\end{equation}
\end{small}%
\endgroup
The steady-state lateral force balance is given as
\par\nobreak\vspace{-5pt}
\begingroup
\allowdisplaybreaks
\begin{small}
\begin{equation}
\label{eqn16}
0 = F_\mathrm{y,R} + F_\mathrm{y,F} - \cos(\phi)\,F_\mathrm{c} + \sin(\phi)\,m\,g,
\end{equation}
\end{small}%
\endgroup
where $F_\mathrm{y,R}$ and $F_\mathrm{y,F}$ are the lateral forces on the rear and front axles, respectively. The vertical force balance is given as
\par\nobreak\vspace{-5pt}
\begingroup
\allowdisplaybreaks
\begin{small}
\begin{multline}
\label{eqn17}
0 = F_\mathrm{z,R} + F_\mathrm{z,F}- \cos(\phi)\cos(\theta)\,m\,g - \sin(\phi)\,\cos(\theta) F_\mathrm{c}\\
 - \underbrace{\frac{1}{m}\,\rho \,C_\mathrm{l}A\,E_\mathrm{kin}}_{F_\mathrm{down}} ,
\end{multline}
\end{small}%
\endgroup
where $F_\mathrm{y,R}$ and $F_\mathrm{y,F}$ are the lateral forces on the rear and front axles, respectively, and $C_\mathrm{l}$ is the downforce coefficient. Given the steady-state cornering assumption, the yaw moment balance is
\par\nobreak\vspace{-5pt}
\begingroup
\allowdisplaybreaks
\begin{small}
\begin{equation}
\label{eqn18}
0 = l_\mathrm{F}\,F_\mathrm{y,F} - l_\mathrm{R}\,F_\mathrm{y,R},
\end{equation}
\end{small}%
\endgroup
where $l_\mathrm{R}$ and $l_\mathrm{F}$ are the front and rear axle longitudinal distances from the center of gravity (CoG). The steady-state pitch moment balance is given as
\par\nobreak\vspace{-5pt}
\begingroup
\allowdisplaybreaks
\begin{small}
\begin{multline}
\label{eqn19}
0 = -l_\mathrm{F}\,F_\mathrm{z,F} + l_\mathrm{R}\,F_\mathrm{z,R} - h_\mathrm{g}
\left(F_\mathrm{x,R}+F_\mathrm{x,F}\right) - l_\mathrm{GP}\frac{1}{m}\,\rho \,C_\mathrm{l}A\,E_\mathrm{kin}\\
- (h_\mathrm{P}-h_\mathrm{G})\frac{1}{m}\,\rho \,C_\mathrm{d}A\,E_\mathrm{kin},
\end{multline}
\end{small}%
\endgroup
where $h_\mathrm{G}$ is the height of the center of gravity, $l_\mathrm{GP}$ is the longitudinal distance between the center of gravity and the center of pressure (CoP) and $h_\mathrm{P}$ is the height of the center of pressure.\\
Tire grip limits are enforced through friction ellipses. These are typically enforced as
\par\nobreak\vspace{-5pt}
\begingroup
\allowdisplaybreaks
\begin{small}
\begin{equation}
\label{eqn20}
F_\mathrm{z,R}\geq \left\|\begin{array}{c}
\frac{F_\mathrm{x,R}}{\mu_\mathrm{x,R}}\\
\frac{F_\mathrm{y,R}}{\mu_\mathrm{y,R}}\\
\end{array}\right\|_2 \quad \text{\normalsize and} \quad F_\mathrm{z,F}\geq \left\|\begin{array}{c}
\frac{F_\mathrm{x,F}}{\mu_\mathrm{x,F}}\\
\frac{F_\mathrm{y,F}}{\mu_\mathrm{y,F}}\\
\end{array}\right\|_2,
\end{equation}
\end{small}%
\endgroup
where $\mu_\mathrm{\{x,y\},\{R,F\}}$ are tire grip coefficients. Often considered constant in previously proposed convex models, these coefficients are, in practice, roughly linearly dependent on load~\cite{Pacejka2012}. Consequently, any load transfer between tires, both longitudinal and lateral, reduces overall grip. To capture this effect, we model the axle friction ellipses as
\par\nobreak\vspace{-5pt}
\begingroup
\allowdisplaybreaks
\begin{small}
\begin{equation}
\label{eqn21}
F_\mathrm{z,R}^\star\geq \left\|\begin{array}{c}
\frac{F_\mathrm{x,R}}{\mu_\mathrm{x,R,nom}}\\
\frac{F_\mathrm{y,R}}{\mu_\mathrm{y,R,nom}}\\
\end{array}\right\|_2 \quad \text{\normalsize and} \quad F_\mathrm{z,F}^\star\geq \left\|\begin{array}{c}
\frac{F_\mathrm{x,F}}{\mu_\mathrm{x,F,nom}}\\
\frac{F_\mathrm{y,F}}{\mu_\mathrm{y,F,nom}}\\
\end{array}\right\|_2,
\end{equation}
\end{small}%
\endgroup
where $F_\mathrm{z,R}^\star$ and $F_\mathrm{z,F}^\star$ are auxiliary variables and $\mu_\mathrm{\{x,y\},\{R,F\},nom}$ are nominal grip coefficients, corresponding to the commonly used Magic Tire Formula parameter $p_\mathrm{D\{x,y\}1}$~\cite{Pacejka2012}. Subsequently, we add the quadratic relations
\par\nobreak\vspace{-5pt}
\begingroup
\allowdisplaybreaks
\begin{small}
\begin{equation}
\label{eqn22}
F_\mathrm{z,R}^\star\leq \frac{\gamma_\mathrm{R}}{2\,F_\mathrm{z,R,nom}}\left(F_\mathrm{z,R}^2+\Delta F_\mathrm{z,R}^2\right) + (1-\gamma_\mathrm{R})F_\mathrm{z,R}
\end{equation}
\end{small}%
\endgroup
and
\par\nobreak\vspace{-5pt}
\begingroup
\allowdisplaybreaks
\begin{small}
\begin{equation}
\label{eqn23}
F_\mathrm{z,F}^\star\leq \frac{\gamma_\mathrm{F}}{2\,F_\mathrm{z,F,nom}}\left(F_\mathrm{z,F}^2+\Delta F_\mathrm{z,F}^2\right) + (1-\gamma_\mathrm{F})F_\mathrm{z,F},
\end{equation}
\end{small}%
\endgroup
\begin{figure}[t]
	\centering
	\includegraphics[width=\linewidth]{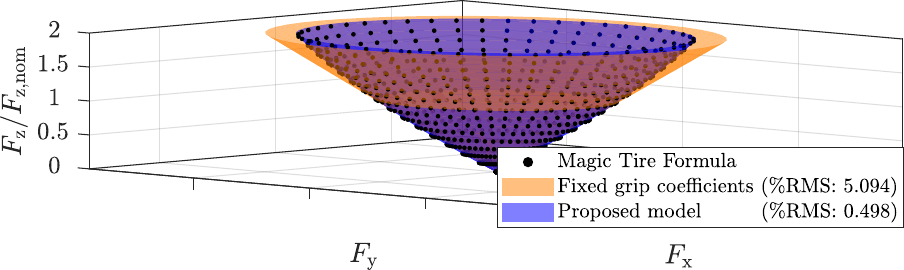}
	\caption{Tire grip envelope of the proposed load-dependent grip coefficients model compared to the simpler model with fixed coefficients and the high-fidelity Magic Formula tire model. Our proposed implementation reduces the normalized root-mean-square error by a factor of 10.}
	\label{fig:gripenvelope}
\end{figure}%
where $\Delta F_\mathrm{z,\{R,F\}}$ is the force difference between the two tires on an axle resulting from lateral load transfer, $F_\mathrm{z,\{R,F\},nom}$ is the vertical force per tire where the nominal grip coefficients apply and $\gamma_\mathrm{\{R,F\}}\leq0$ describes the decrease in grip coefficients over increased load. This roughly follows from Magic Tire Formula parameters with $\gamma \approx \frac{p_\mathrm{Dy2}}{p_\mathrm{Dy1}}$~\cite{Pacejka2012}. This inequality holds with equality when the grip limit constraint is active, as $F_\mathrm{z,\{R,F\}}^\star$ is maximized. The force differences over the tires on an axle are defined as
\par\nobreak\vspace{-5pt}
\begingroup
\allowdisplaybreaks
\begin{small}
\begin{equation}
\label{eqn24}
\Delta F_\mathrm{z,R} = \frac{2\,h_\mathrm{rc,R}\,F_\mathrm{y,R}}{w_\mathrm{R}}+(1-\xi)\frac{2\left(h_\mathrm{G}F_\mathrm{y}-h_\mathrm{rc,R}F_\mathrm{y,R}-h_\mathrm{rc,F}F_\mathrm{y,F}\right)}{w_\mathrm{R}}
\end{equation}
\end{small}%
\endgroup
and
\par
\begin{figure}[t]
	\centering
	\includegraphics[width=0.6\linewidth]{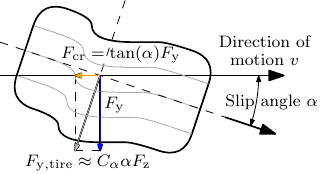}
	\caption{Origins of the tire cornering resistance. No lateral force is possible without a tire slip angle, and the slip angle creates a force component opposing the direction of motion.}
	\label{fig:slipangle}
\end{figure}

\nobreak\vspace{-5pt}
\begingroup
\allowdisplaybreaks
\begin{small}
\begin{equation}
\label{eqn25}
\Delta F_\mathrm{z,F} = \frac{2\,h_\mathrm{rc,F}\,F_\mathrm{y,F}}{w_\mathrm{F}}+\xi\,\frac{2\left(h_\mathrm{G}F_\mathrm{y}-h_\mathrm{rc,R}F_\mathrm{y,R}-h_\mathrm{rc,F}F_\mathrm{y,F}\right)}{w_\mathrm{F}},
\end{equation}
\end{small}%
\endgroup
where $w_\mathrm{\{R.F\}}$ is the track width on the respective axle, $h_\mathrm{rc,\{R,F\}}$ is the roll center height, and $\xi$ is the roll stiffness distribution. These equations follow from the steady-state roll moment balance. The resulting feasible tire force space is visualized in Fig.~\ref{fig:gripenvelope}, along with the simpler model with fixed tire grip coefficients and the corresponding Magic Tire Formula operating space. The proposed tire model matches much better with the high-fidelity Magic Formula model, reducing the normalized root-mean-square error by a factor of 10, while retaining the second-order conic form.
Another common omission in convex lap-time optimization is a cornering-induced resistance force. The origins of this effect are visualized in Fig.~\ref{fig:slipangle}. To develop a lateral force, the tire needs a \textit{slip angle} $\alpha$. This induces a force component opposing the direction of motion that is proportional to the tangent of the slip angle. Assuming small slip angles, we can linearize the tangent, and the slip angle is proportional to $\frac{F_\mathrm{y}}{F_\mathrm{z}}$~\cite{Pacejka2012}. Thereby, we model the resistance induced by cornering as 
\par\nobreak\vspace{-5pt}
\begingroup
\allowdisplaybreaks
\begin{small}
\begin{equation}
\label{eqn26}
F_\mathrm{cr,\{R,F\}} = \frac{F_\mathrm{y,\{R,F\}}^2}{C_\mathrm{\alpha,\{R,F\}}F_\mathrm{z,\{R,F\}}},
\end{equation}
\end{small}%
\endgroup
where $F_\mathrm{cr,\{R,F\}}$ is the cornering resistance and $C_\mathrm{\alpha,\{R,F\}}\leq0$ is the cornering stiffness, equivalent to $p_\mathrm{Ky1}$ in the Magic Formula definitions~\cite{Pacejka2012}. Although this is naturally nonconvex, we can subtract it from powertrain constraints to ensure convexity. Herein we also include the rolling resistances as follows:
\par\nobreak\vspace{-5pt}
\begingroup
\allowdisplaybreaks
\begin{small}
\begin{equation}
\label{eqn27}
F_\mathrm{x,\{R,F\}} \leq F_\mathrm{w,\{R,F\}} - c_\mathrm{r,\{R,F\}}F_\mathrm{z,\{R,F\}}+ \frac{F_\mathrm{y,\{R,F\}}^2}{C_\mathrm{\alpha,\{R,F\}}F_\mathrm{z,\{R,F\}}},
\end{equation}
\end{small}%
\endgroup
where $F_\mathrm{w,\{R,F\}}$ is the force exerted on the wheel by the brakes and powertrain, and $c_\mathrm{r,\{R,F\}}$ is the rolling resistance coefficient per axle. This constraint is only convex for $F_\mathrm{z,\{R,F\}}\geq0$. In order to make it suitable for SOCP solvers, we rewrite it in the convex second-order conic form as
\par\nobreak\vspace{-5pt}
\begingroup
\allowdisplaybreaks
\begin{small}
\begin{multline}
\label{eqn28}
-F_\mathrm{x,\{R,F\}} +F_\mathrm{w,\{R,F\}} - \left(c_\mathrm{r,\{R,F\}}+C_\mathrm{\alpha,\{R,F\}}\right)F_\mathrm{z,\{R,F\}} \geq \\ \left\|\begin{array}{c}
2\,F_\mathrm{y,\{R,F\}}\\
F_\mathrm{x,\{R,F\}} - F_\mathrm{w,\{R,F\}} + \left(c_\mathrm{r,\{R,F\}}-C_\mathrm{\alpha,\{R,F\}}\right)F_\mathrm{z,\{R,F\}}\\
\end{array}\right\|_2.
\end{multline}
\end{small}%
\endgroup
Finally, torque and power upper limitations are enforced as
\par\nobreak\vspace{-5pt}
\begingroup
\allowdisplaybreaks
\begin{small}
\begin{equation}
\label{eqn29}
F_\mathrm{w,\{R,F\}} \leq P_\mathrm{max,\{R,F\}}{\textstyle\frac{\mathrm{d}t}{\mathrm{d}s}}
\end{equation}
\end{small}%
\endgroup
and
\par\nobreak\vspace{-5pt}
\begingroup
\allowdisplaybreaks
\begin{small}
\begin{equation}
\label{eqn30}
F_\mathrm{w,\{R,F\}} \leq  \frac{T_\mathrm{max,\{R,F\}}}{r_\mathrm{w,\{R,F\}}},
\end{equation}
\end{small}%
\endgroup
where $T_\mathrm{max,\{R,F\}}$ is the maximum torque per axle, $r_\mathrm{w,\{R,F\}}$ is the wheel radius per axle, and $P_\mathrm{max,\{R,F\}}$ is the maximum propulsive power per axle. Lower limits may be imposed in the same way.
\subsection{Optimization Problem}
We formulate the control problem using the state variables $\mathbf{x}(s_\mathrm{ref})  \coloneq \{n,n',E_\mathrm{kin}\}$ and the control variables $\mathbf{u}(s_\mathrm{ref})  \coloneq 
\{n'',F_\mathrm{w,R},F_\mathrm{w,F}\}$. Herein $n''$ is related to the steering input and the wheel forces follow from brake and powertrain operation.
\begin{prob}[Minimum-time Optimization Problem]\label{Problem1}
The minimum-lap-time control strategy is the solution of\\
\\
$\begin{aligned}
\min_{\mathbf{u}(s_\mathrm{ref}),\mathbf{x}(s_\mathrm{ref})} \quad 
  &t_{\mathrm{lap}}\\[0.75em]
\text{subject to } \quad &\\
    \begin{rcases}
        (\ref{eqn1}), (\ref{eqn3}), (\ref{eqn4}), (\ref{eqn5}-\ref{eqn8}), \\
        (\ref{eqn14}-\ref{eqn19}), (\ref{eqn24}), (\ref{eqn25}), \\
        \mathbf{x}(0) = \mathbf{x}(S_\mathrm{ref})
      \end{rcases} &h(\mathbf{x},\mathbf{u}) = 0 \\[1em]
    \begin{rcases}
        (\ref{eqn2}), (\ref{eqn11}), (\ref{eqn12}), (\ref{eqn13}),\\
        (\ref{eqn21}-\ref{eqn23}), (\ref{eqn28}-\ref{eqn30}) \\
      \end{rcases} &g(\mathbf{x},\mathbf{u}) \leq 0.
\end{aligned}$
\\
\end{prob}%
\noindent Several equations here remain nonconvex, namely (\ref{eqn5}), (\ref{eqn7nl}), (\ref{eqn9}) and (\ref{eqn15nl}). This means we cannot take advantage of specialized convex programming methods for computational efficiency. Therefore, we propose an approach similar to Heilmeier et al.~\cite{Heilmeier2019} to rewrite these as affine constraints, by taking a first-order Taylor approximation using results from a previous sequential programming iteration. We linearize the centrifugal force equation (\ref{eqn5}) to
\par\nobreak\vspace{-5pt}
\begingroup
\allowdisplaybreaks
\begin{small}
	\begin{equation}
	\label{eqn6}
	F_\mathrm{c} = 2E_{\mathrm{kin},k-1}\left.\left[
	\begin{matrix}
	\nabla_{x'} \kappa\\
	\nabla_{y'} \kappa\\
	\nabla_{x''} \kappa\\
	\nabla_{y''} \kappa\\
	\end{matrix}
	\right]\right|_{k-1}^\top
	\left[
	\begin{matrix}
x'_k-x'_{k-1}\\
y'_k-y'_{k-1}\\
x''_k-x''_{k-1}\\
y''_k-y''_{k-1}\\
	\end{matrix}	
	\right]+2E_{\mathrm{kin},k}\kappa_{k-1},
	\end{equation}
\end{small}%
\endgroup
where $k$ is the iteration. These approximations work well since the racing line curvature does not deviate significantly from the reference curvature of the circuit center, which itself provides a baseline initial guess for the first iteration. We base the circuit slope (\ref{eqn7nl}) on our previous iteration result,
\par\nobreak\vspace{-5pt}
\begingroup
\allowdisplaybreaks
\begin{small}
	\begin{equation}
		\label{eqn7}
		\theta(s_\mathrm{ref}) = \arctan{\left(\frac{z'_{k-1}(s_\mathrm{ref})}{\sqrt{x'^2_{k-1}(s_\mathrm{ref})+y'^2_{k-1}(s_\mathrm{ref})}}\right)},
	\end{equation}
\end{small}%
\endgroup
and we linearize our bilinear domain transformations (\ref{eqn9}) and (\ref{eqn15nl}):
\par\nobreak\vspace{-5pt}
\begingroup
\allowdisplaybreaks
\begin{small}
	\begin{equation}
		\label{eqn10}
	t_{\mathrm{lap},k} = \int_{0}^{S_\mathrm{ref}}{\smash{\textstyle\frac{\mathrm{d}t}{\mathrm{d}s}}}_k{\smash{\textstyle\frac{\mathrm{d}s}{\mathrm{d}s_\mathrm{ref}}}}_{k-1}+{\smash{\textstyle\frac{\mathrm{d}t}{\mathrm{d}s}}}_{k-1}{\smash{\textstyle\frac{\mathrm{d}s}{\mathrm{d}s_\mathrm{ref}}}}_k-{\smash{\textstyle\frac{\mathrm{d}t}{\mathrm{d}s}}}_{k-1}{\smash{\textstyle\frac{\mathrm{d}s}{\mathrm{d}s_\mathrm{ref}}}}_{k-1}\, \mathrm{d}s_\mathrm{ref}
	\end{equation}
\end{small}%
\endgroup
and
\par\nobreak\vspace{-5pt}
\begingroup
\allowdisplaybreaks
\begin{small}
	\begin{multline}
		\label{eqn15}
		{\smash{\textstyle\frac{\mathrm{d}E_\mathrm{kin}}{\mathrm{d}s_\mathrm{ref}}}} = {\smash{\textstyle\frac{\mathrm{d}E_\mathrm{kin}}{\mathrm{d}s}}}_k{\smash{\textstyle\frac{\mathrm{d}s}{\mathrm{d}s_\mathrm{ref}}}}_{k-1}+{\smash{\textstyle\frac{\mathrm{d}E_\mathrm{kin}}{\mathrm{d}s}}}_{k-1}{\smash{\textstyle\frac{\mathrm{d}s}{\mathrm{d}s_\mathrm{ref}}}}_k-{\smash{\textstyle\frac{\mathrm{d}E_\mathrm{kin}}{\mathrm{d}s}}}_{k-1}{\smash{\textstyle\frac{\mathrm{d}s}{\mathrm{d}s_\mathrm{ref}}}}_{k-1}.
	\end{multline}
\end{small}%
\endgroup
\renewcommand{\algorithmicendwhile}{\textbf{end}}
\begin{algorithm}[t]
	\caption{Iterative Solving Procedure for Problem~\ref{Problem1}(SCP)}\label{alg1}
	\begin{algorithmic}
		%\STATE \textbf{inputs:}\\
		\STATE $\mathbf{x}_0,\mathbf{u}_0 \gets \text{Initial Guess}$ \COMMENT{can be constant over $s_\mathrm{ref}$} \\
		\STATE $k = 1$
		\WHILE{$\|\mathbf{x}_{k}-\mathbf{x}_{k-1}\| > \epsilon \text{ (threshold)}$} 
		\STATE $\mathbf{x}_k,\mathbf{u}_k \gets  \text{Solve Problem~\ref{Problem2} linearized around $\mathbf{x}_{k-1},\mathbf{u}_{k-1}$}$\\
		\STATE $k = k + 1$
		\ENDWHILE
	\end{algorithmic}
\end{algorithm}%
Subsequently, our convexified, iteratively solved subproblem reads as follows:
\begin{prob}[Convex Subproblem]\label{Problem2}
	The convexified minimum-lap-time control strategy at iteration $k$ is the solution of\\
	\\
$\begin{aligned}
	\min_{\mathbf{u}_k(s_\mathrm{ref}),\mathbf{x}_k(s_\mathrm{ref})} \quad 
	&t_{\mathrm{lap},k}\\[0.75em]
	\text{subject to } \quad &\\
	\begin{rcases}
		(\ref{eqn1}), (\ref{eqn3}), (\ref{eqn4}),(\ref{eqn8}),(\ref{eqn14}),(\ref{eqn16}-\ref{eqn19}),\\
		(\ref{eqn24}),(\ref{eqn25}),(\ref{eqn6}-\ref{eqn15}) \\
		\mathbf{x}_k(0) = \mathbf{x}_k(S_\mathrm{ref})
	\end{rcases} &h_k(\mathbf{x},\mathbf{u}) = 0 \\[1em]
	\begin{rcases}
		(\ref{eqn2}), (\ref{eqn11}), (\ref{eqn12}), (\ref{eqn13}),(\ref{eqn21}-\ref{eqn23}), (\ref{eqn28}-\ref{eqn30}) \\
	\end{rcases} &g_k(\mathbf{x},\mathbf{u}) \leq 0.
\end{aligned}$
	\\
\end{prob}%
\noindent Since this subproblem consists exclusively of linear equalities and second-order conic inequalities (or a subset thereof), we can use dedicated SOCP solvers to exploit this model structure and reduce solving time~\cite{Domahidi2013}. The iterative sequential-convex programming procedure to solve Problem~\ref{Problem1} is outlined in Algorithm~\ref{alg1}. The problem is considered solved to a satisfactory level once the norm of the state change is less than a predetermined tolerance $\epsilon$.

\subsection{Remarks on Convergence}
Due to its convexity, Problem~\ref{Problem2} is consistently solved with global optimality guarantees and requires no informed initialization. However, with the presented implementation, convergence of the iterative procedure outlined in Algorithm~\ref{alg1} is not guaranteed.
If the linearized variables deviate significantly from their previous iterate linearization point, the system could destabilize. A typical solution in sequential-convex programming that does not change the model structure is to add an adaptive \textit{trust region} constraint that limits this deviation, guaranteeing convergence~\cite{Martins_Ning_2021}. In practice, for this particular problem, this was found to be unnecessary in all tested cases, even with an uninformed initial guess.
The narrow track boundary constraint \eqref{eqn2} inherently limits  offsets of linearized variables similar to how a trust region would (which is also why the linearization method is effective in the first place). For unmonitored applications, a general trust region mechanism within the loop of Algorithm~\ref{alg1} is advisable, though this may come at the cost of requiring additional solving iterations. Since such methods are well-established in literature and our contribution centers on the proposed modeling framework, we consider their detailed treatment beyond the scope of this work.

\begin{figure}[t]
	\centering
	\includegraphics[width=\linewidth]{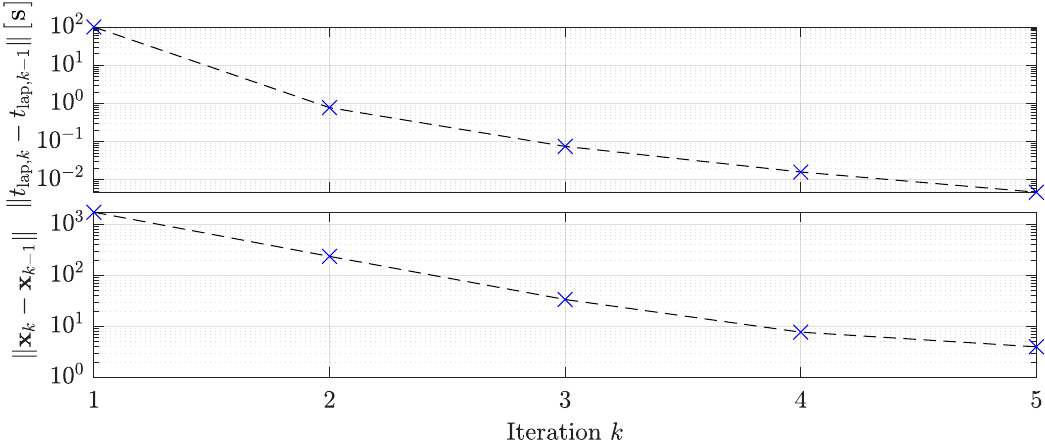}
	\caption{Convergence of the SCP iterative procedure. Without a warm-start, only 4 iterations are needed to converge within one-hundredth of a second.}
	\label{fig:convergence}
\end{figure}
\section{Results}
\label{ResultsSection}
This section presents numerical results regarding the performance
of our proposed sequential convex programming method. We base our results on a Formula 1 race car driving on the Spa-Francorchamps circuit. From regulations and other publicly available sources, we can make realistic estimates of car parameters~\cite{FIA2025TechRegs}. Circuit data is acquired from rFactor2~\cite{rFactor2}, a high-fidelity racing simulator. For the discretization, we apply the Hermite-Simpson collocation method with 2000 samples, resulting in a fixed step size of around 3.5m. We parse the problem using YALMIP~\cite{Lofberg2004} in MATLAB, and solve the convex problem using ECOS~\cite{Domahidi2013}. We solve the problem on a standard consumer laptop running Windows with an Intel Core i7-11800H 2.3 GHz processor. Thereby, we achieve consistent solving times of around 1 second per iteration. We deem the problem solved to a satisfactory level when the lap-time improvement is less than one-hundred of a second. Generally, this occurs after 4 to 5 iterations. Thereby, even without a warm start, the full problem is solved in less than 5 seconds. In fact, as shown in Fig.~\ref{fig:convergence}, it converges at a rate where just a single iteration could suffice in a model-predictive-control (MPC) setting where the solution is close to an initial guess, paving the way for real-time trajectory planning applications. Moreover, mesh refinement techniques could further reduce the number of required samples and thereby computation time.
\begin{figure}[t]
	\centering
	\includegraphics[width=\linewidth]{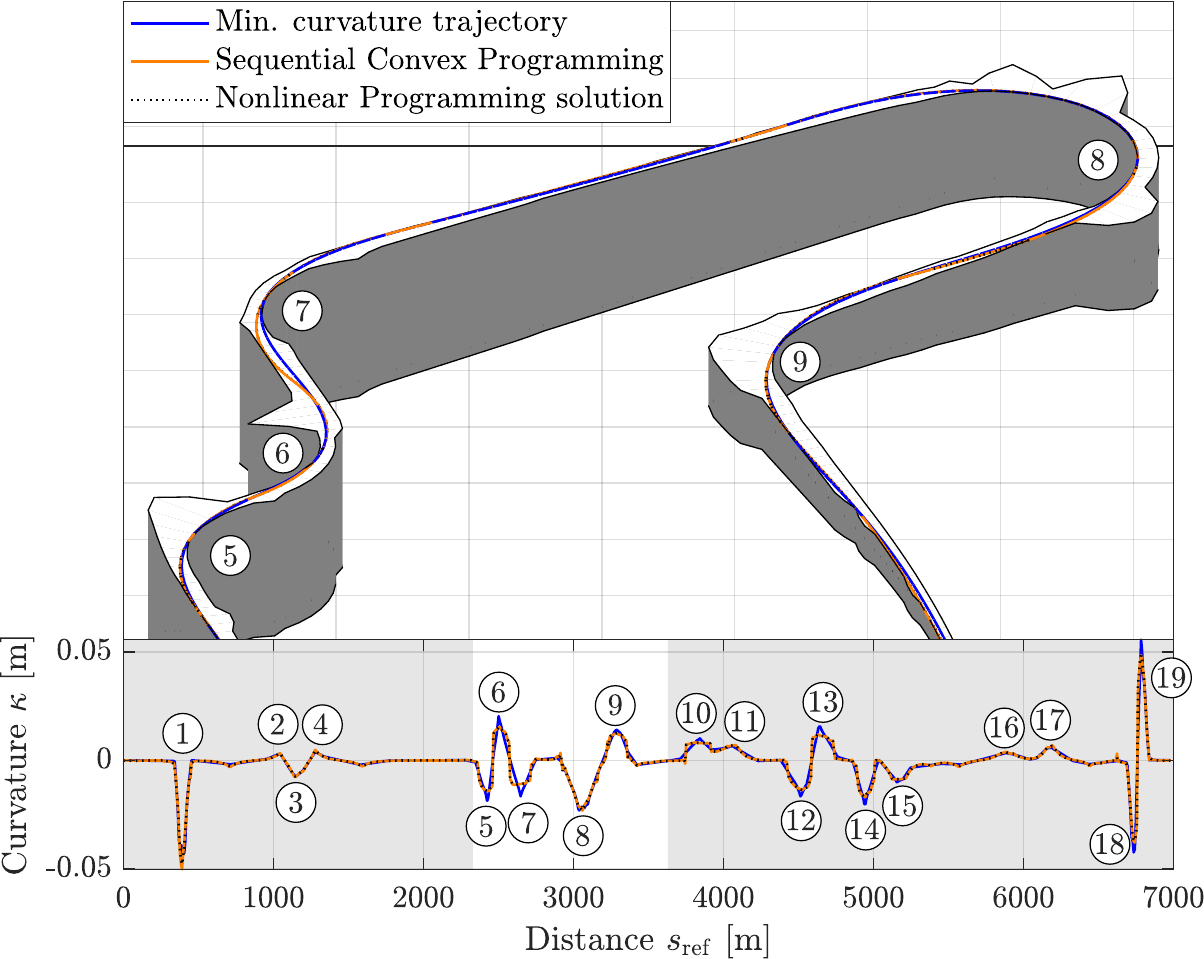}
	\caption{Racing line visualization through turns 5 to 9 on the Spa-Francorchamps circuit, comparing the fast SCP and minimum-curvature solutions to the computationally heavy NLP solution. The SCP solution matches the NLP solution, and highlights the weakness of the fixed-trajectory minimum-curvature approach.}
	\label{fig:trajectory}
\end{figure}
\subsection{Trajectory Results}
\label{TrajectoryResultsSection}
\begin{figure}[t]
	\centering
	\includegraphics[width=\linewidth]{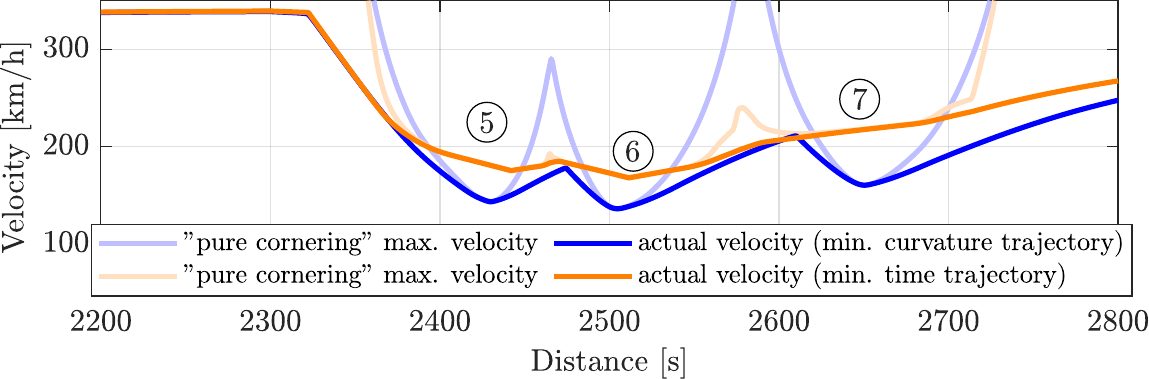}
	\caption{The minimum-curvature trajectory yields a faster "pure cornering" speed profile (meaning all tire grip is devoted to maximizing instantaneous cornering speed). This potential is not exploited due to acceleration limitations, making the minimum-time trajectory significantly quicker.}
	\label{fig:curvatureexploitation}
\end{figure}

\begin{figure*}[t]
	\centering
	\includegraphics[width=\linewidth]{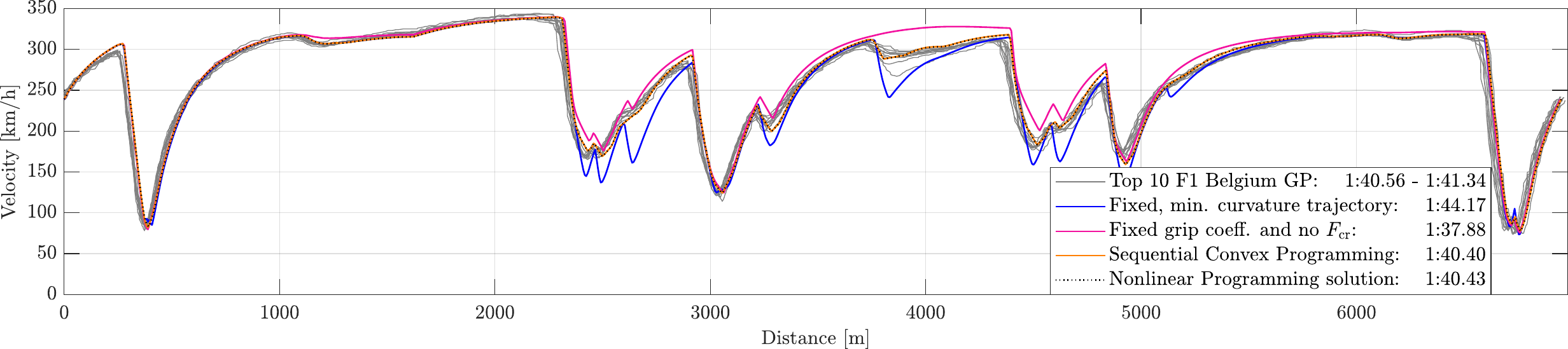}
	\caption{The SCP solution matches the NLP solution. Both solutions match the real-world data well. The minimum curvature racing line method's limitations are apparent, as its total lap-time is almost 4 seconds off the free trajectory lap-time. The model with fixed grip coefficients and no cornering resistance overestimates grip at high speed and fails to capture deceleration in some high-speed curved sections, highlighting the importance of these additions.}
	\label{fig:velocities}
\end{figure*}

The resulting spacial trajectory is visualized in Fig.~\ref{fig:trajectory}. The result of our fast sequential-convex programming method solving Problem~\ref{Problem2} matches that of Problem~\ref{Problem1} directly solved through nonlinear programming with IPOPT~\cite{Waechter2005}. However, the nonlinear programming method takes more than 2 minutes to solve the problem despite using the same tolerances. This means that the SCP method has a lower solving time of around a factor of 25, demonstrating the significant computation time benefits. This figure also highlights the shortcomings of the fixed-trajectory approach based on a preprocessed minimum-curvature racing line~\cite{Heilmeier2019,Kapania2019}. This method achieves efficient solving times by ignoring vehicle limitations and separately optimizes the racing line on a minimum-curvature basis, which operates on the rationale that a minimized curvature maximizes speed potential. We see that the minimum-curvature solution prefers short peaks in curvature, whereas the free-trajectory minimum-lap-time-based solutions prefer a longer sustained curvature through the corners. The peaks in curvature determine the minimum speed through the corners, and the minimum-curvature solution is unable to fully exploit the quick drop-off in curvature (i.e., velocity potential) due to acceleration and braking limitations (see Fig.~\ref{fig:curvatureexploitation}). As shown in Fig.~\ref{fig:velocities}, the lower cornering speeds of the minimum-curvature trajectory result in a significant final lap-time deficit of around 4\% compared to minimum-lap-time-based trajectory optimizations. This figure again shows the match between the fast sequential \
convex programming solution and the computationally heavy nonlinear model solution. In addition, a result is included that highlights the importance of cornering-induced resistances and load-dependent grip coefficients, which are often omitted from convex models. The proposed method that considers these effects is better capable of estimating grip at higher speeds and captures details like reduced acceleration at full power when driving a slightly curved path. For a real-world reference, we also compare our model to the fastest laps of the top 10 qualifiers in the 2025 Formula One Belgian GP~\cite{Schaefer2025FastF1}. To avoid overfitting to this dataset, we only tuned the drag- and downforce coefficient to get this level of correlation (well within a reportedly realistic range). Thereby, even with limited correlation work, a decent match is achieved to real-world data, with our optimization result deviating no more than 10 km/h.
\subsection{g-g-v Diagram}
\begin{figure*}[t]
	\centering
	\includegraphics[width=\linewidth]{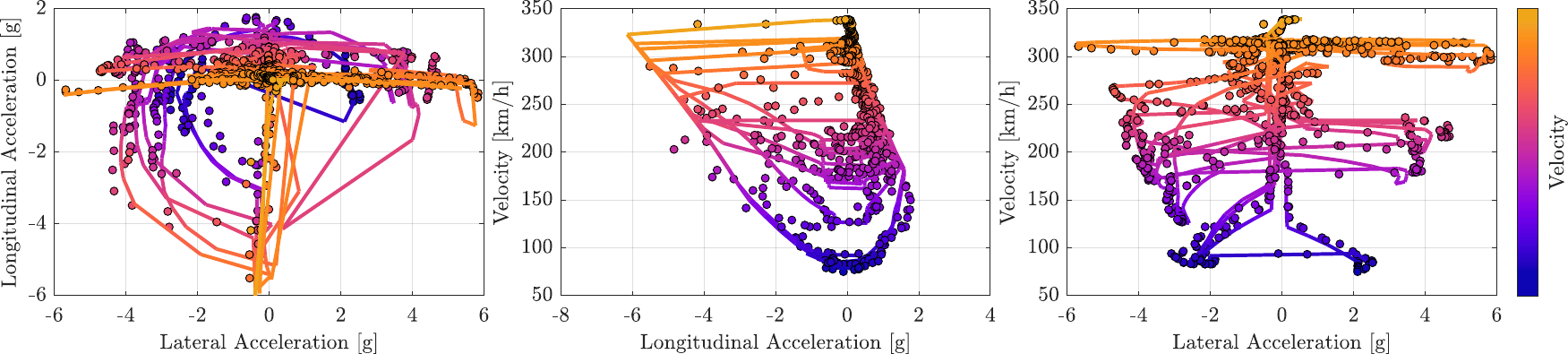}
	\caption{g-g-v diagram: the optimization model solution (line trace) compared to real-world data (markers).}
	\label{fig:ggv}
\end{figure*}

As mentioned in the introduction, a typical method of performing efficient lap-time optimization is by precomputing a so-called g-g-v diagram~\cite{Lovato2021,Lovato2021a,Duhr2022}. Such a diagram holds information on where the vehicle-dynamic performance boundaries lie~\cite{Milliken1995}. As such, for an optimal lap-time, one would expect the car to be at the boundary throughout the lap. Fig.~\ref{fig:ggv} shows such a diagram for the SCP optimization result compared with data from the fastest Formula One lap, assuming both drive the same racing line. The g-g-v diagram produced by the optimization is of the typical shape and captures the details well. This demonstrates that the vehicle-dynamic model outlined in Section~\ref{ModelSection} parametrically defines the g-g-v boundaries and is representative of real-world limits, validating the model potential.
\subsection{Energy Management}
\begin{figure}[t]
	\centering
	\includegraphics[width=\linewidth]{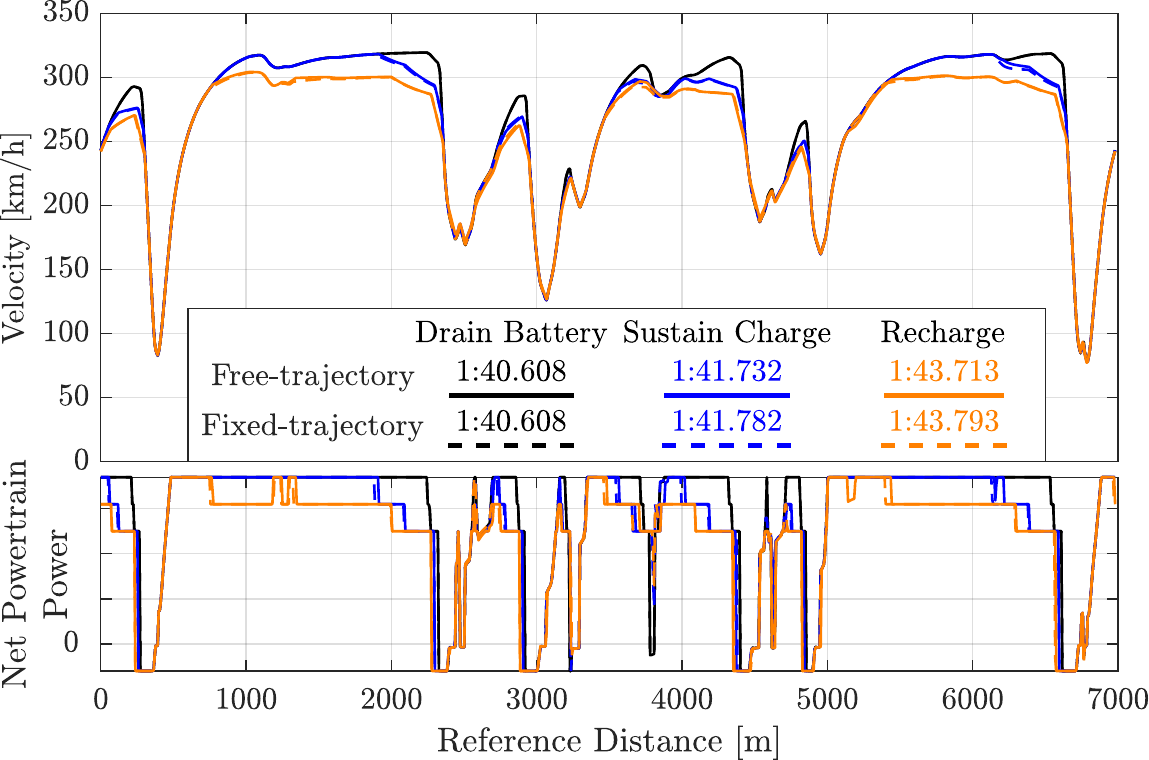}
	\caption{Three different energy management strategies on a hybrid-electric racing car, comparing optimizations with a fixed trajectory based on the battery drain scenario to the respective free-trajectory solution. There is a minor difference in lap-time between fixed- and free-trajectory solutions, showing that the driven trajectory is somewhat influenced by the energy management strategy.}
	\label{fig:energymanagement}
\end{figure}

In addition to the improved trajectory generation method, our framework is also suitable to determine optimal powertrain operation and energy management. An additional state may be introduced as the distance integral of drivetrain forces such as $F_\mathrm{w}$, which captures the cumulative energy flow. Imposing limitations on this energy will result in an optimized energy management solution. Hybrid powertrain topologies are possible by linearly splitting $F_\mathrm{w}$ into, for example, braking forces, electric motor forces and combustion engine forces, each with their own limitations~\cite{DuhrChristodoulouEtAl2020,EbbesenSalazarEtAl2018,Kampen2024}. We demonstrate this in Fig.~\ref{fig:energymanagement}, where we simultaneously investigate the commonly used assumption that an energy management strategy is not significantly affected by the chosen racing line and vice versa~\cite{DuhrChristodoulouEtAl2020,EbbesenSalazarEtAl2018,Kampen2024,EshofKampenEtAl2025}. We perform six optimizations on a hybrid racing car. In addition to its combustion engine, it features a capacity-limited battery with an electric motor capable of delivering a limited extra power boost and regenerative braking. We consider three scenarios with different energy constraints: One starts with a full battery and may fully drain it during the lap, one starts with an empty battery and has to fill it by the end of the lap, and one has to sustain its battery charge, matching it at the beginning and end. Each scenario is evaluated with a free trajectory and with a fixed trajectory based on the battery drain scenario. Naturally, the results are identical for the battery drain scenario, as it is using the same trajectory. For the charge-sustaining and battery-filling scenario, there is a marginal difference in lap-time of under 0.1\%, quantifying that energy limitations do influence the driving trajectory to a minor extent.
\section{Conclusion} \label{Conclusion}
This paper presented a sequential convex programming approach for free-trajectory lap-time optimization.
Our method is capable of jointly optimizing the racing line and powertrain operation orders of magnitude faster than standard nonlinear programming algorithms. This makes it suitable for global path planning in autonomous racing, design studies and parameter sensitivity analysis.
We demonstrated how equally computationally efficient minimum-curvature methods can lead to significantly slower lap-times compared to our proposed minimum-time optimization, and demonstated our framework's ability to optimize energy-limited powertrain operation.
Finally, we validated the common assumption in energy management studies that energy limitations do not alter the racing line and found it results in only minor lap-time losses.

In future work, we would like to leverage our framework for real-time trajectory planning and energy management purposes.

\section*{Acknowledgment}
\noindent We thank Dr. I. New for proofreading this paper.

\bibliographystyle{IEEEtran}
\renewcommand{\baselinestretch}{0.96}
\bibliography{../../../bibliography/powertrains,../../../bibliography/main,../../../bibliography/SML_papers}

\end{document}